\documentclass[reqno,12pt]{amsart}
\theoremstyle{plain}    
\newtheorem{thm}{Theorem}[section]
\newtheorem{cor}[thm]{Corollary} 
\newtheorem{lemma}[thm]{Lemma} 
\newtheorem{prop}[thm]{Proposition}

\theoremstyle{remark}

\theoremstyle{definition}
\newtheorem{conj}[thm]{Conjecture}

\usepackage{amscd,amssymb,comment,eepic,euscript,graphics}

\newcommand\alphat{{\tilde{\alpha}}}
\newcommand{\bb}[1]{{\mathbb{#1}}}
\newcommand\Bc{{\mathcal{B}}}
\newcommand\bt{{\tilde b}}
\newcommand{\cl}[1]{{\mathcal{#1}}}
\newcommand\Cpx{\bb{C}}
\newcommand\ct{{\tilde c}}
\newcommand\eps{\epsilon}
\newcommand\HEu{{\EuScript H}}                   

\newcommand\Ints{{\bb Z}}
\newcommand\Lt{{\widetilde{L}}}
\newcommand\Mcal{\cl M}
\newcommand\Nats{{\bb N}}
\newcommand\oneh{{\hat 1}}
\newcommand\ran{\operatorname{ran}}

\newcommand\St{\widetilde{S}}
\newcommand\Tt{{\widetilde{T}}}
\newcommand\xh{{\hat x}}
\newcommand\cc{^\ast}

\def\d{{\rm d}}

\def\R{\mathbb{R}}

\begin{document}

\title[Brown measures and Aluthge iterates for some operators]
{Brown measure and iterates of the Aluthge transform for some operators arising from measurable actions$^*$}

\author[Dykema, Schultz]
{Ken Dykema$^\flat$, Hanne Schultz$^{\dag\,\ddag}$}

\address{K. Dykema, Department of Mathematics, Texas A\&M University,
College Station, TX 77843-3368, USA}
\email{kdykema@math.tamu.edu}

\address{H. Schultz,
Department of Mathematics and Computer Science,
University of Southern Denmark,
Campusvej 55,
5230 Odense M, Denmark}
\email{schultz@imada.sdu.dk}

\thanks{$^*$An earlier version of this paper was distributed under the title:
On the Aluthge transform: continuity properties and Brown measure.}

\thanks{$^\flat$Research supported in part by NSF grant DMS--0300336.}

\thanks{$^\dag$As a student of the Ph.D.\ school OP-ALG-TOP-GEO this author is partially
supported by the Danish Research Training Council.}
\thanks{$^\ddag$Partially supported by The Danish National Research Foundation.}

\subjclass[2000]{47A05, (47B99)}

\keywords{Aluthge transform, Brown measure, mean ergodic theorem}

\date{\today}

\begin{abstract}
We consider the Aluthge transform $\Tt=|T|^{1/2}U|T|^{1/2}$ of a
Hilbert space operator $T$,
where $T=U|T|$ is the polar decomposition of $T$.
We prove that the map $T\mapsto\Tt$ is continuous with respect to the norm
topology and with respect to the $*$--SOT topology on bounded sets.
We consider the special case in a tracial von Neumann algebra when
$U$ implements
an automorphism of the von Neumann algebra generated by the positive part $|T|$ of $T$, and
we prove that the iterated Aluthge transform converges to a normal operator whose
Brown measure agrees with that of $T$ (and we compute this Brown measure).
This proof relies on a theorem that is
an analogue of von Neumann's mean ergodic theorem, but for sums weighted by
binomial coefficients.
\end{abstract}

\maketitle

\section{Introduction}

For a Hilbert space $\HEu$, we let $\Bc(\HEu)$ denote the set of bounded linear
operators on $\HEu$.
Let $T\in\Bc(\HEu)$ and let $T=U|T|$ be its polar decomposition.
The {\em Aluthge transform} of $T$ is the operator $\Tt=|T|^{1/2}U|T|^{1/2}$.
This was first studied in~\cite{Al90} and has received much attention in recent years;
(see, for example, the papers cited below).
One reason the Aluthge transform is interesting is in relation to the invariant subspace
problem.
Jung, Ko and Pearcy prove in~\cite{JKP00} that $T$ has a nontrivial invariant subspace
if and only if $\Tt$ does.
They also note that when $T$ is a quasiaffinity, then $T$ has a nontrivial hyperinvariant
subspace if and only if $\Tt$ does.
(A quasiaffinity is a operator with zero kernel and dense range;
the invariant and hyperinvariant subspace problems are interesting only for quasiaffinities.)
Clearly, the spectrum of $\Tt$ equals that of $T$.
Jung, Ko and Pearcy prove in~\cite{JKP00}
that other spectral data are also preserved by the Aluthge transform;
see also~\cite{JKP01} for more and related results.

The iterated Aluthge transforms (or Aluthge iterates)
of $T$ are the operators $\Tt^{(k)}$, $k\ge0$,
defined by setting $\Tt^{(0)}=T$ and letting $\Tt^{(k+1)}$ be the Aluthge transform
of $\Tt^{(k)}$.
An interesting result, due to Yamazaki~\cite{Y02} (see also~\cite{W03}),
is that the spectral radius of $T$ is equal to the limit as $k\to\infty$
of the norm of $\Tt^{(k)}$.
Recently,
It was conjectured in~\cite{JKP03} that $\Tt^{(k)}$ converges in strong operator topology
to a normal operator as $k\to\infty$.
Antezana, Pujals and Stojanoff~\cite{APS07} proved this conjecture for operators on finite
dimensional Hilbert spaces.
While this conjecture has been shown in~\cite{CJL} not to hold in general, we believe,
based partly on the evidence of examples described below, that it will
hold for elements of finite von Neumann algebras.

The examples we consider are when $T=U|T|$, where $U$ is unitary and where
the positive part $|T|$ belongs to an abelian algebra that is normalized by $U$,
(i.e.\ on which $U$ acts by conjugation).
In this situation, we describe the Aluthge iterates of $T$ and show that the conjecture
mentioned above holds.
We also caculate the Brown measure for these operators.
Our proof relies on an analogue of von Neumann's mean ergodic theorem for sums
weighted with binomial coefficients, and on a theorem of Haagerup and Schultz.

The authors thank the referee for identifying some
needed improvements in the paper.

\section{Continuity properties in $\Bc(\HEu)$}
\label{sec:cont}

The results of this section are applications of the continuous functional calculus
for positive operators.

\begin{lemma}\label{lem:fnT}
For $n\in\Nats$ consider the function
$f_n(t)=\sqrt{\max(\frac1n,t)}$ defined for real, non--negative $t$
and let $A_n=f_n(|T|)$.
Then
\begin{align}
\|A_n\|&\le\max\big({\textstyle \frac1{\sqrt n}},\|T\|^{1/2}\big) \label{eq:fnT1}\\
\|\,|T|A_n^{-1}\|&\le\|T\|^{1/2} \\
\|A_n-|T|^{1/2}\|&\le{\textstyle \frac1{\sqrt n}} \\
\|\,|T|A_n^{-1}-|T|^{1/2}\|&\le{\textstyle \frac1{4\sqrt n}}  \label{eq:fnT4} \\
\|A_nTA_n^{-1}-\Tt\|&\le{\textstyle \frac5{4\sqrt n}}\|T\|^{1/2}. \label{eq:fnT5}
\end{align}
\end{lemma}
\begin{proof}
The norm bounds~\eqref{eq:fnT1}--\eqref{eq:fnT4} follow easily from the
continuous functional calculus.
For~\eqref{eq:fnT5}, we use
\[
\|A_nTA_n^{-1}-\Tt\|\le
\begin{aligned}[t]
&\|A_n-|T|^{1/2}\|\,\|\,U|T|A_n^{-1}\| \\
&+\|\,|T|^{1/2}U\|\,\|\,|T|A_n^{-1}-|T|^{1/2}\|
\end{aligned}
\]
and apply the previous estimates.
\end{proof}

\begin{lemma}\label{lem:pq}
Given $R\ge1$ and $\eps>0$, there are real polynomials $p$ and $q$
such that for every $T\in\Bc(\HEu)$ with $\|T\|\le R$, we have
\[
\|\Tt-p(T^*T)Tq(T^*T)\|<\eps.
\]
\end{lemma}
\begin{proof}
Take $\eps<1$, for convenience.
By Lemma~\ref{lem:fnT}, we may choose $n$ so large that
\[
\|\Tt-A_nTA_n^{-1}\|<\eps
\]
whenever $\|T\|\le R$, where $A_n=f_n(|T|)$, with $f_n$ as in Lemma~\ref{lem:fnT}.
Let $p$ and $q$ be real polynomials of one variable such that
\[
\max_{t\in[0,R]}|f_n(t)-p(t^2)|<\eps,\qquad\max_{t\in[0,R]}|f_n(t)^{-1}-q(t^2)|<\eps.
\]
Then we have
\[
\|A_n-p(T^*T)\|<\eps,\qquad\|A_n^{-1}-q(T^*T)\|<\eps.
\]
Using the estimates in Lemma~\ref{lem:fnT}, we get
\begin{multline*}
\|A_nTA_n^{-1}-p(T^*T)Tq(T^*T)\| \\
\le\|A_n-p(T^*T)\|\,\|TA_n^{-1}\|+\|p(T^*T)T\|\,\|A_n^{-1}-q(T^*T)\|<3\eps R^{3/2}.
\end{multline*}
Now replace $\eps$ in the above argument by $\eps/(3R^{3/2})$.
\end{proof}

In the statements below, $\|\cdot\|$ refers to the operator norm topology and
$*$--SOT to the $*$--strong operator topology on $\Bc(\HEu)$.

\begin{thm}\label{thm:cont}
The Aluthge transform map $T\mapsto\Tt$ is:
\begin{itemize}
\item[(a)]
$(\|\cdot\|,\|\cdot\|)$--continuous on $\Bc(\HEu)$,
\item[(b)] ($*$--SOT,$*$--SOT)--continuous on bounded subsets of $\Bc(\HEu)$.
\end{itemize}
\end{thm}
\begin{proof}
For~(a), let $T\in\Bc(\HEu)$ and take $\eps\in(0,1]$.
Let $R=\|T\|+1$ and let $p$ and $q$ be polynomials as found in Lemma~\ref{lem:pq}
for these values of $R$ and $\eps$.
Let $\delta\in(0,1]$ be such that $\|T-S\|<\delta$ implies
\[
\|p(T^*T)Tq(T^*T)-p(S^*S)Sq(S^*S)\|<\eps.
\]
Then $\|T-S\|<\delta$ implies
\[
\|\Tt-\St\|<\|\Tt-p(T^*T)Tq(T^*T)\|+\|\St-p(S^*S)Sq(S^*S)\|+\eps<3\eps,
\]
where the last inequality is by choice of $p$ and $q$.
This proves part~(a).

For~(b), let $R>1$, $\eps>0$
and let $p$ and $q$ be the polynomials found in Lemma~\ref{lem:pq} for these values.
Let $x\in\HEu$.
If $S,T\in\Bc(\HEu)$, each with norm $\le R$, then
\begin{align*}
\|(\Tt-\St)x\|&\le2\eps\|x\|+\|(p(T^*T)Tq(T^*T)-p(S^*S)Sq(S^*S))x\|, \\
\|(\Tt-\St)^*x\|&\le2\eps\|x\|+\|(p(T^*T)T^*q(T^*T)-p(S^*S)S^*q(S^*S))x\|,
\end{align*}
Since multiplication is $*$--SOT continuous on bounded subsets of $\Bc(\HEu)$,
by taking $S$ in some neighborhood of $T$ in this topology, both quantities
$\|(p(T^*T)Tq(T^*T)-p(S^*S)Sq(S^*S))x\|$
and
$\|(p(T^*T)T^*q(T^*T)-p(S^*S)S^*q(S^*S))x\|$
can be forced to be arbitrarily small.
This proves (b).
\end{proof}

\section{The Aluthge transform in finite von Neumann algebras}
\label{sec:fvN}

In this section, we consider a von Neumann algebra $\Mcal$ equipped with
a normal, faithful, tracial state $\tau$, acting on the Hilbert space
$\HEu:=L^2(\Mcal,\tau)$,
which is the completion of $\Mcal$ with respect to the norm $\|x\|_2=\tau(x^*x)^{1/2}$.
For $x\in\Mcal$, we denote the corresponding element of $L^2(\Mcal,\tau)$ by $\xh$.
Clearly, the Aluthge transform $\Tt$ of any $T\in\Mcal$ also lies in $\Mcal$.
In $\Mcal$, convergence in SOT in $B(\HEu)$
implies convergence in $\|\cdot\|_2$, because $\|x\|_2=\|x\oneh\|_{\HEu}$.
On the other hand, on bounded subsets of $\Mcal$, convergence in $\|\cdot\|_2$ implies convergence in $*$--SOT.
Indeed, let $a_n$ be a sequence in the unit ball of $\Mcal$ such that $\|a_n\|_2\to0$
as $n\to\infty$,
let $\eta\in L^2(\Mcal,\tau)$, and let us show $\|a_n\eta\|_\HEu\to0$ as $n\to\infty$.
Let $\eps>0$.
Then there is $x\in\Mcal$ such that $\|\eta-\xh\|_\HEu<\eps$.
Letting $\rho$ denote the right action of $\Mcal^{op}$ on $\HEu$, we have
\[
\|a_n\eta\|_\HEu\le2\eps+\|a_n\xh\|_\HEu=2\eps+\|\rho(x)a_n\oneh\|\HEu\le2\eps+\|x\|\|a_n\|_2\,,
\]
and from this we conclude $\|a_n\eta\|_\HEu\to0$ as $n\to\infty$.
Since $\|a_n^*\|_2=\|a_2\|_2$, we conclude $a_n\to0$ in $*$--SOT.

Therefore, we have the following immediate corollary of Theorem~\ref{thm:cont}(b).
\begin{thm}\label{thm:2nmconts}
The Aluthge transformation $T\mapsto\Tt$ is $(\|\cdot\|_2,\|\cdot\|_2)$--continuous
on bounded subsets of $\Mcal$.
\end{thm}

\vspace{.2cm}

The proof of Theorem~2.2 of~\cite{JKP03} carries over to the setting
of finite von Neumann algebras to prove the following result.
\begin{thm}\label{thm:2-norm}
If $T\in\Mcal$, then
\begin{equation}\label{eq:2-norm}
\|\Tt\|_2\le\|T\|_2.
\end{equation}
Moreover, equality holds if and only if $T$ is normal.
\end{thm}

\vspace{.2cm}

\begin{prop}\label{prop:ifConv}
Let $T\in\Mcal$.
If any subsequence $\{\Tt^{(n_k)}\}_{k=1}^\infty$ of the sequence of Aluthge interates of $T$
converges in $\|\cdot\|_2$ to a limit $L$,
then $L$ is a normal operator.
\end{prop}
\begin{proof}
By~\eqref{eq:2-norm}, $\|\Tt^{(n)}\|_2$ decreases to $\|L\|_2$ as $n\to\infty$.
By Theorem~\ref{thm:2nmconts}, $\Tt^{(n_k+1)}$ converges in $\|\cdot\|_2$
as $k\to\infty$ to $\Lt$.
Therefore, $\|\Lt\|_2=\|L\|_2$ and, by Theorem~\ref{thm:2-norm}, $L$ is normal.
\end{proof}

\vspace{.2cm}

\begin{lemma}\label{lem:An-1}
Let $T\in\Mcal$ and let $A_n=f_n(|T|)$ be as defined in Lemma~\ref{lem:fnT}.
Let $P_0\in\Mcal$ be the projection onto $\ker T$.
Then
\begin{equation}\label{eq:1P0}
\lim_{n\to\infty}\|A_n^{-1}|T|^{1/2}-(1-P_0)\|_2=0.
\end{equation}
\end{lemma}
\begin{proof}
We have $A_n^{-1}|T|^{1/2}=h_n(|T|)$, where
\[
h_n(t)=
\begin{cases}
\sqrt nt^{1/2},&t\in[0,\frac1n] \\
1,&t\in[\frac1n,\infty).
\end{cases}
\]
Therefore, $0\le(1-P_0)-A_n^{-1}|T|^{1/2}\le E_{|T|}((0,\frac1n))$,
where $E_{|T|}$ is the spectral measure of $|T|$.
Since $\big(E_{|T|}((0,\frac1n))\big)_{n=1}^\infty$ is a sequence of projections decreasing to zero as $n\to\infty$,
we get~\eqref{eq:1P0}.
\end{proof}


The Brown measure~\cite{B86} of an operator $T\in\Mcal$ is defined as $\frac1{2\pi}$ times the Laplacian
of the function $\Cpx\ni\lambda\mapsto\Delta(\lambda-T)$, where $\Delta$ is the Fuglede--Kadison
determinant~\cite{FK52}, defined by $\Delta(X)=\exp(\tau(\log|X|))$.
By Theorem~4.3 of~\cite{B86}, the Brown measures $\mu_T$ and $\mu_{\Tt}$ are the same.
In light of this, one would like in Proposition~\ref{prop:ifConv}
to have that $L$ has the same Brown measure as $T$.
This is at present unknown.
However, we will conjecture even more than this,
namely the actual convergence of Aluthge iterates.

\begin{conj}\label{bigconj}
Let $T\in \Mcal$.
Then the sequence $\{\Tt^{(n)}\}_{n=1}^\infty$ of Aluthge iterates of $T$
converges in $\|\cdot\|_2$ to a normal operator $N$ whose Brown measure is equal
to the Brown measure of $T$. 
\end{conj}

This conjecture is the analogue for finite von Neumann algebras of Conjecture~5.6 of~\cite{JKP03},
which is about SOT--convergence of Aluthge iterates of $T\in\Bc(\HEu)$.
Although Cho, Jung and Lee~\cite{CJL} solve Conjecture~5.6 of~\cite{JKP03}
in the negative by giving an example of a weighted shift operator $T$
whose Aluthge iterates fail to converge even in weak operator topology,
that example clearly generates an infinite von Neumann algebra.
As remarked earlier, this conjecture has recently been proved for matrices in~\cite{APS07}.

\section{An Ergodic theorem for sums weighted with binomial coefficients}
\label{sec:ergodic}

In this section, we prove a result that resembles
von Neumann's mean ergodic theorem, but for sums
that are weighted with binomial coefficients.
Our proof is based on the proof of von Neumann's theorem found in~\cite{P83}.
This theorem is covered in the case of a measure--preserving transformation
by Theorem~1 of~\cite{HP}.
However, check of the relevant conditions in our case seems to be just as much work as proving the theorem directly.
After proving the theorem, we also draw some consequences that will be used in some examples in Section~\ref{sec:ex}.

\begin{thm}\label{thm:meanergodic}
Let $T\in\Bc(\HEu)$ satisfy $\|T\|\le1$.
Then for every $v\in\HEu$,
the vector
\[
H_n(v)=\frac1{2^n}\sum_{k=0}^n\binom nkT^kv
\]
converges as $n\to\infty$ to $Pv$, where $P$ is the orthogonal
projection of $\HEu$ onto its subspace $\ker(T-I)=\{x\in\HEu\mid Tx=x\}$.
\end{thm}
\begin{proof}
We have
\begin{equation}\label{eq:ker}
\ker(T-I)=\ker(T^*-I).
\end{equation}
Indeed, it will suffice to show the inclusion $\subseteq$.
If $Tx=x$, then
\begin{align*}
0\le\|T^*x-x\|^2&=\|T^*x\|^2-\langle T^*x,x\rangle-\langle x,T^*x\rangle+\|x\|^2 \\
&=\|T^*x\|^2-\langle x,Tx\rangle-\langle Tx,x\rangle+\|x\|^2
=\|T^*x\|^2-\|x\|^2\le0,
\end{align*}
and~\eqref{eq:ker} is proved.
Therefore, $\ker(T-I)=\ran(T-I)^\perp$.

It is clear that if $v\in\ker(T-I)$, then $H_n(v)=v$ for all $n$.
We will show that if $v\in\ker(T-I)^\perp$, then $\|H_n(v)\|$ converges to zero
as $n\to\infty$.
By linearity, this will imply that for general $v\in\HEu$, $H_n(v)$ converges to $Pv$,
In order to show that $\|H_n(v)\|$ converges to zero for $v\in\ker(T-I)^\perp$,
it will suffice to show it for $v\in\ran(T-I)$, because this latter space is
dense in $\ker(T-I)^\perp$ and each $H_n$ is a linear contraction.
So we may assume $v=Ty-y$, some $y\in\HEu$.
Then
\begin{align*}
H_n(v)&=\frac1{2^n}\sum_{k=0}^n\binom nk(T^{k+1}y-T^ky) \\
&=-\frac1{2^n}y+\frac1{2^n}\sum_{k=1}^n\bigg(\binom n{k-1}-\binom nk\bigg)T^ky
+\frac1{2^n}T^{n+1}y.
\end{align*}
Since
\[
\bigg|\binom n{k-1}-\binom nk\bigg|=\frac{|2k-n-1|}k\binom n{k-1},
\]
we get
\[
\|H_n(v)\|\le\frac1{2^n}\bigg(2+\sum_{k=1}^n\frac{|2k-n-1|}k\binom n{k-1}\bigg)\|y\|.
\]
It will suffice to show
\begin{equation}\label{eq:binomsum}
\lim_{n\to\infty}\frac1{2^n}\sum_{k=1}^n\frac{|2k-n+1|}k\binom nk=0.
\end{equation}
Let $0<\alpha<\frac12$.
Using Stirling's formula, it is easy to prove that for some constant $C>0$, we always have
\begin{equation}
\sum_{k=n-\lfloor\alpha n-1\rfloor}^n\binom nk=\sum_{k=0}^{\lfloor\alpha n-1\rfloor}\binom nk
\le C\,n\bigg(\frac1{\alpha^\alpha(1-\alpha)^{(1-\alpha)}}\bigg)^n.
\end{equation}
(And better estimates are possible; see, for example,~\cite{E83}.)
Since
\[
\alpha^\alpha(1-\alpha)^{(1-\alpha)}>\frac12,
\]
we have
\begin{equation}\label{eq:limsup}
\limsup_{n\to\infty}\frac1{2^n}\sum_{k=1}^n\frac{|2k-n+1|}k\binom nk
=\limsup_{n\to\infty}\frac1{2^n}
 \sum_{k=\lfloor\alpha n\rfloor}^{n-\lfloor\alpha n\rfloor}\frac{|2k-n+1|}k\binom nk.
\end{equation}
However, for $\alpha n-1\le k\le n-\alpha n+1$, we have
\[
\frac{|2k-n+1|}k\le\frac{(1-2\alpha)n+3}{\alpha n-1},
\]
which together with~\eqref{eq:limsup} yields
\[
\limsup_{n\to\infty}\frac1{2^n}\sum_{k=1}^n\frac{|2k-n+1|}k\binom nk\le\frac{1-2\alpha}\alpha.
\]
Taking $\alpha$ arbitrarily close to $\frac12$ gives~\eqref{eq:binomsum}.
\end{proof}

\vspace{.2cm}

Let us now consider a probability space $(X,\mu)$
and a $\mu$--preserving, invertible transformation $\alpha:X\to X$.
Let $E^\alpha:L^1(\mu)\to L^1(\mu)^\alpha$ denote the conditional expectation onto the
subspace $L^1(\mu)^\alpha$ of $\alpha$--invariant functions in $L^1(\mu)$.
In analogous notation, we will also write
$E^\alpha:L^2(\mu)\to L^2(\mu)^\alpha$ and $E^\alpha:L^\infty(\mu)\to L^\infty(\mu)^\alpha$
for the restrictions of $E^\alpha$ to the indicated subspaces.
For a random variable $b$ and $n\ge0$, we will consider the random variable
\begin{equation}\label{eq:hn}
h_n(b)=\frac1{2^n}\sum_{k=0}^n\binom nkb\circ\alpha^k.
\end{equation}
Applying Theorem~\ref{thm:meanergodic} in the case $\HEu=L^2(\mu)$ and $Tf=f\circ\alpha$,
we have the following result.
\begin{cor}\label{cor:2ndmom}
Let $b\in L^2(\mu)$ and let $c=E^\alpha(b)$.
Then $\lim_{n\to\infty}\|h_n(b)-c\|_2=0$.
\end{cor}

\vspace{.2cm}

\begin{cor}\label{cor:1stmom}
Let $b\in L^1(X)$ and let $c=E^\alpha(b)$.
Then $\lim_{n\to\infty}\|h_n(b)-c\|_1=0$.
\end{cor}
\begin{proof}
For $M>0$, let $g_M:[-\infty,\infty]\to[-M,M]$ be
\begin{equation}\label{eq:gM}
g_M(t)=
\begin{cases}
-M,&t\le-M \\
t,&-M\le t\le M \\
M,&M\le t.
\end{cases}
\end{equation}
Let $\eps>0$.
For $M>0$ sufficiently large, letting $\bt=g_M\circ b$
we have $\|\bt-b\|_1<\eps$, so $\|h_n(\bt)-h_n(b)\|_1<\eps$.
Moreover, we have $\|\ct-c\|_1<\eps$, where $\ct=E^\alpha(\bt)$.
Since we're working over a probability space,
from Corollary~\ref{cor:2ndmom} we have
\[
\lim_{n\to\infty}\|h_n(\bt)-\ct\|_1=0,
\]
and this yields
\[
\limsup_{n\to\infty}\|h_n(b)-c\|_1\le2\eps.
\]
\end{proof}

We now consider a random variable $b:X\to[-\infty,R]$ for some real number $R$.
Let $L^\infty(X)^\alpha$ denote the space of functions $f\in L^\infty(X)$ satisfying $f\circ\alpha=f$.
The conditional expectation of $b$, $E^\alpha(b):X\to[-\infty,R]$, is the unique $\alpha$--invariant random variable
such that $\int E^\alpha(b)f\,d\mu=\int bf\,d\mu\in[-\infty,R]$ whenever $f\in L^\infty(X)^\alpha$.
Then the map $b\mapsto E^\alpha(b)$ is linear and order preserving.
Moreover, $b$ is integrable if and only if $E^\alpha(b)$ is integrable.

\begin{cor}\label{cor:notL1}
Suppose $R$ is a real number and $b:X\to[-\infty,R]$ is measurable.
Then $h_n(b)$ converges in probability to $E^\alpha(b)$ as $n\to\infty$,
where we take the obvious metric on $[-\infty,R]$.
\end{cor}
\begin{proof}
Fix $K>0$ large and $\eps>0$ and $\delta>0$ small.
Consider the $\alpha$--invariant set
$F=\{\omega\mid E^\alpha(b)(\omega)\ge-K\}$.
Then $E^\alpha(1_Fb)=1_FE^\alpha(b)$ is integrable, so $1_Fb$ is integrable.
By Corollary~\ref{cor:1stmom}, $\lim_{n\to\infty}\|h_n(1_Fb)-E^\alpha(1_Fb)\|_1=0$.
Also, we have $h_n(1_Fb)=1_Fh_n(b)$.
Thus, for $n$ sufficiently large we have
\begin{equation}\label{eq:delta}
\mu(\{\omega\mid E^\alpha(b)(\omega)\ge-K,\,|h_n(b)(\omega)-E^\alpha(b)(\omega)|\ge\eps\})<\delta.
\end{equation}

For integers $m$ with $m>|R|$ and with $g_m$ as defined in~\eqref{eq:gM},
we have $E^\alpha(g_m\circ b)\ge E^\alpha(g_{m+1}\circ b)\ge E^\alpha(b)$
and we claim that $E^\alpha(g_m\circ b)$ converges
a.e.\ as $m\to\infty$ to $E^\alpha(b)$.
Indeed, let
$H=\lim_{m\to\infty}E^\alpha(g_m\circ b)$.
Choosing $f\in L^\infty(X)^\alpha$ nonnegative, we have by monotone convergence
\[
\int E^\alpha(b)f\,d\mu=\int b\,f\,d\mu=\lim_{m\to\infty}\int(g_m\circ b)\,fd\mu
=\lim_{m\to\infty}\int E^\alpha(g_m\circ b)\,f\,d\mu=\int H\,f\,d\mu.
\]
This implies $H=E^\alpha(b)$ a.e.

By Egoroff's Theorem, $\exp(E^\alpha(g_m\circ b))$ converges in measure to $\exp(E^\alpha(b))$,
so $E^\alpha(g_m\circ b)$ converges in measure to $E^\alpha(b)$.
Therefore, for some $m$ we have
\[
\mu(\{\omega\mid E^\alpha(b)(\omega)<-K,\,E^\alpha(g_m\circ b)(\omega)\ge-K+1\})<\delta.
\]
Again by Corollary~\ref{cor:1stmom}, $\lim_{n\to\infty}\|h_n(g_m\circ b)-E^\alpha(g_m\circ b)\|_1=0$,
so for $n$ sufficiently large we have
\[
\mu(\{\omega\mid E^\alpha(b)(\omega)<-K,\,h_n(g_m\circ b)\ge-K+2\})<2\delta.
\]
Finally, since $h_n(b)\le h_n(g_m\circ b)$, for $n$ sufficiently large we have
\begin{equation}\label{eq:2delta}
\mu(\{\omega\mid E^\alpha(b)(\omega)<-K,\,h_n(b)\ge-K+2\})<2\delta.
\end{equation}
Combining~\eqref{eq:delta} and~\eqref{eq:2delta} finishes the proof.
\end{proof}

Corollaries~\ref{cor:1stmom} and~\ref{cor:notL1} are, of course, straightforward
consequences of Theorem~\ref{thm:meanergodic}.
Analogues where von Neumann's mean ergodic theorem is used instead of Theorem~\ref{thm:meanergodic}
are obtained similarly.
For future use, we state the following such analogue of Corollary~\ref{cor:notL1}.

\begin{prop}\label{prop:meanergodic}
Suppose $R$ is a real number and $b:X\to[-\infty,R]$ is measurable.
Then the random variable
\[
\frac1n\sum_{k=0}^{n-1}b\circ\alpha^k
\]
converges in probability to $E^\alpha(b)$ as $n\to\infty$.
\end{prop}

\section{Some examples in finite von Neumann algebras}
\label{sec:ex}

Let $T\in\Bc(\HEu)$ and let $T=U|T|$ be its polar decomposition.
The Aluthge transform of $T$ is by definition $\Tt=|T|^{1/2}U|T|^{1/2}$.
It is easily seen that if $V\in\Bc(\HEu)$ is any partial isometry whose restriction to
$U^*U\HEu$ equals $U$, then $\Tt = |T|^{1/2}V|T|^{1/2}$.
In particular, if $T$ belongs to a finite von Neumann algebra, then we may take $V$ to be
a unitary extension of $U$.

We consider a probability space $(X,\mu)$
and we let $\alpha$ be an invertible, measure--preserving transformation of $X$.
Let $\Mcal=L^\infty(X)\rtimes_\alpha\Ints$ be the crossed product von Neumann algebra,\
which is, therefore, generated by a copy of $L^\infty(X)=L^\infty(X,\mu)$
and a unitary $U$ such that for $f\in L^\infty(X)$, we have
$UfU\cc = f\circ\alpha$.
For convenience, we will write
\[
\alphat(f)=f\circ\alpha^{-1}=U^*fU,
\]
so that we have $fU=U\alphat(f)$.
The set 
\[
{\rm span}\{U^k f\,|\,k\in\Ints, f\in L^\infty(X)\}
\]
is then strongly dense in $\Mcal$, and there is a normal tracial state $\tau$ on $\Mcal$ uniquely determined by
\begin{equation}\label{trace}
\tau(U^k f)=\delta_{k,0}\int_X f\d\mu, \qquad k\in\Ints, f\in L^\infty(X).
\end{equation}
We will investigate the Aluthge iterates of operators of the form $T=U|T|\in\Mcal$ with $U$ as  above and with $|T|\in L^\infty(X)$.

\begin{lemma}\label{absvalue}
The $n$'th Aluthge iterate of $T$ is
\begin{equation}\label{eq:TtU}
\Tt^{(n)}=U|\Tt^{(n)}|,
\end{equation}
where
\begin{equation}\label{eq1}
|\Tt^{(n)}|=\prod_{k=0}^n\alphat^{k}(|T|^{\binom nk/2^n}).
\end{equation}
\end{lemma}
\begin{proof} The proof proceeds by induction over $n$.
Clearly,~\eqref{eq:TtU} and~\eqref{eq1} hold for $n=0$
(with the convention that $\Tt^{(0)}=T$).
Let $N\geq 1$, and assume that \eqref{eq1} holds for $n=N-1$.
Then for the $N$th Aluthge iterate we have:
\begin{align*}
\Tt^{(N)} &= \Big[\prod_{k=0}^{N-1}\alphat^{k}(|T|^{\binom {N-1}k/2^{N-1}})\Big]^{\frac12} U  \Big[\prod_{l=0}^{N-1}\alphat^{l}(|T|^{\binom {N-1}l/2^{N-1}})\Big]^{\frac12}\\
&=  \Big[\prod_{k=0}^{N-1}\alphat^{k}(|T|^{\binom {N-1}k/2^{N}})\Big] U  \Big[\prod_{l=0}^{N-1}\alphat^{l}(|T|^{\binom {N-1}l/2^{N}})\Big] \displaybreak[2] \\
&= U \Big[\prod_{k=0}^{N-1}\alphat^{k+1}(|T|^{\binom {N-1}k/2^{N}})\Big]\Big[\prod_{l=0}^{N-1}\alphat^{l}(|T|^{\binom {N-1}l/2^{N}})\Big] \displaybreak[2] \\
&= U \Big[\prod_{k=1}^{N}\alphat^{k}(|T|^{\binom {N-1}{k-1}/2^{N}})\Big]\Big[\prod_{l=0}^{N-1}\alphat^{l}(|T|^{\binom {N-1}l/2^{N}})\Big] \displaybreak[2] \\
&= U \alphat^{N}(|T|^{\binom {N-1}{N-1}/2^{N}}) \Big[\prod_{k=1}^{N-1}\alphat^{k}(|T|^{\binom {N-1}{k-1}/2^N+\binom {N-1}k/2^{N}})\Big]|T|^{\frac{1}{2^N}}\\
&= U \prod_{k=0}^N \alphat^{k}(|T|^{\binom Nk/2^N}),
\end{align*}
and this shows that~\eqref{eq:TtU} and~\eqref{eq1} hold for $n=N$ as well. 
\end{proof}

Given a random variable $b:X\to[-\infty,R]$ for some real number $R$,
we let $E^\alpha(b)$ denote the conditional expectation of $b$
onto the space of such $\alpha$--invariant random variables, as described
immediately before Corollary~\ref{cor:notL1}.

\begin{thm}\label{thm:AluthgeConv}
Let $T=U|T|$ be as described above.
Then
\[ \lim_{n\to\infty}\|\Tt^{(n)}-UH\|_2=0, \]
where $H=\exp(E^\alpha(\log(|T|)))$.
\end{thm}
\begin{proof} In light of~\eqref{eq:TtU},
it will suffice to show $|\Tt^{(n)}|\stackrel{\|\cdot\|_2}{\longrightarrow}H$ as $n\to\infty$.
Since the Aluthge iterates $(\Tt^{(n)})_{n=1}^\infty$ form a norm--bounded sequence in~ $L^\infty(X)$,
it clearly suffices to show that $|\Tt^{(n)}|$ converges in probability to $H$.
This is equivalent to convergence in probability of
$\log|\Tt^{(n)}|$ to $E^\alpha(\log(|T|))$.
From~\eqref{eq1},
\[
\log|\Tt^{(n)}|=\frac{1}{2^n}\sum_{k=0}^n \binom nk \log(|T|)\circ\alpha^{k}=h_n(\log(|T|)),
\]
with $h_n$ as defined in~\eqref{eq:hn}.
Now the required convergence follows from Corollary~\ref{cor:notL1}.
\end{proof}

\vspace{.2cm}

If $\alpha$ is ergodic, then $E^\alpha(b)=\tau(b)$ and we have the following.
\begin{cor}
If $T=U|T|$ as above and if $\alpha$ is ergodic, then the Aluthge iterates $(\Tt^{(n)})_{n=1}^\infty$
converge in $\|\cdot\|_2$ to $\Delta(T)U$, where $\Delta(T)\in[0,\infty)$
is the Fuglede--Kadison determinant of $T$.
\end{cor}

\vspace{.2cm}

We will now verify Conjecture~\ref{bigconj} for our operators $T$
by proving that
the Brown measure of $T$ is the same as that of $UH$, with $H$ as in Theorem~\ref{thm:AluthgeConv}.
This is an application of very general results from~\cite{HS2}.

\begin{thm}\label{Brown measure} 
Let $T=U|T|$ and $H=\exp(E^\alpha(\log|T|))$ be as in Theorem~\ref{thm:AluthgeConv}.
Then in the strong operator topology, 
\begin{equation}\label{convtoH 1}
\big[(T^m)\cc T^m\big]^{\frac{1}{2m}}\stackrel{m\rightarrow\infty}{\longrightarrow}H,
\end{equation}
and 
\begin{equation}\label{convtoH 2}
\big[T^m (T^m)\cc\big]^{\frac{1}{2m}}\stackrel{m\rightarrow\infty}{\longrightarrow}H.
\end{equation}
Moreover, $UH$ is a normal operator and
the Brown measures of $T$ and $UH$ agree.
\end{thm}

\begin{proof}
We will prove \eqref{convtoH 1}.
Note that
\[
T^m = U^m\prod_{k=0}^{m-1}\alphat^{k}(|T|),
\]
so that
\begin{equation}\label{eq2}
[(T^m)\cc T^m]^\frac{1}{2m}= \Big[\prod_{k=0}^{m-1}\alphat^{k}(|T|)\Big]^{\frac1m}.
\end{equation}
and
\[
\log[(T^m)\cc T^m]^\frac{1}{2m}=\frac1m\sum_{k=0}^{m-1}\log(\alphat^{k}(|T|))
\]
By Proposition~\ref{prop:meanergodic}, the random variable $\log[(T^m)\cc T^m]^\frac{1}{2m}$
converges in probability to $E^\alpha(\log|T|)$, which implies convergence in probability of
$[(T^m)\cc T^m]^\frac{1}{2m}$ to $H$.
Since the sequence of random variables $([(T^m)\cc T^m]^\frac{1}{2m})_{m\in\Nats}$
is uniformly bounded, convergence in probability implies converges in $\|\cdot\|_2$ and, hence,
also in SOT to $H$.
The SOT--convergence \eqref{convtoH 2} follows analogously.

Since $H$ is $\alpha$-invariant, $H$ commutes with $U$ and $UH$ is normal.
We will now show that the Brown measures $\mu_T$ and $\mu_{UH}$ agree.
First note that both measures are invariant under rotations. Indeed, for $\theta\in\R$ we have an isomorphism $\pi_\theta:\Mcal\rightarrow\Mcal$ given by
\begin{eqnarray*}
\pi_\theta(U)&=& e^{i\theta}U,\\
\pi_\theta(f)&=& f, \qquad f\in L^\infty(X).
\end{eqnarray*}
It is also clear from \eqref{trace} that $\pi_\theta$ preserves the trace so that
\[
\mu_T=\mu_{\pi(T)}=\mu_{e^{i\theta}T} \quad {\rm and} \quad \mu_{UH}=\mu_{\pi(UH)}=\mu_{e^{i\theta}UH}.
\]
Hence, it suffices to show that the two measures agree on the closed disks centered at the origin,
$\overline{B(0,r)}$, for $r>0$.
Since $UH$ is normal,
\[
\mu_{UH}(\overline{B(0,r)})=\tau(1_{[0,r]}(H)).
\]
On the other hand, applying
Thm.~8.1 of~\cite{HS2},
the convergence~\eqref{convtoH 1} and Lemma~7.16(iii) of~\cite{HS2}, we get
\[
\mu_T(\overline{B(0,r)})= \tau(1_{[0,r]}(H)).
\]
\end{proof}

\bibliographystyle{plain}

\begin{thebibliography}{19}

\bibitem{Al90} A.\ Aluthge,
`On $p$--hyponormal operators for $0<p<1$,'
{\em Integr. Equ. Oper. Theory} {\bf 13} (1990), 307-315.

\bibitem{APS07} J.\ Antezana, E.R.\ Pujals, D.\ Stojanoff,
`The iterated Aluthge transforms of a matrix converge,'
preprint, arXiv:0711.3727.

\bibitem{B86} L.G.\ Brown,
`Lidskii's Theorem in the type II case,'
{\em Geometric methods in operator algebras (Kyoto, 1983)},
H.\ Araki and E.G.\ Effros, (Eds.),
Pitman Res. Notes Math. Ser.\ {\bf 123}, Longman Sci.\ Tech., 1986, pp.\ 1-35.

\bibitem{CJL} M.\ Cho, I.\ Jung and W. Lee,
`On Aluthge transforms of $p$--hyponormal operators,'
{\em Integr. Equ. Oper. Theory}, {\bf 53} (2005), 321-329.

\bibitem{E83} J.\ Elton,
`Sign--embeddings of $\ell^n_1$,'
\emph{Trans.\ Amer.\ Math.\ Soc.}\ {\bf 279} (1983), 113-124.

\bibitem{FK52} B.\ Fuglede and R.V.\ Kadison,
`Determinant theory in finite factors,'
{\em Ann.\ of Math.\ (2)} {\bf 55} (1952) 520-530.

\bibitem{JKP00} I.\ Jung, E.\ Ko and C.\ Pearcy,
`Aluthge transforms of operators,'
{\em Integr. Equ. Oper. Theory} {\bf 37} (2000), 437-448.

\bibitem{JKP01} I.\ Jung, E.\ Ko and C.\ Pearcy,
`Spectral pictures of Aluthge transforms of operators,'
{\em Integr. Equ. Oper. Theory} {\bf 40} (2001), 52-60.

\bibitem{JKP03} I.\ Jung, E.\ Ko and C.\ Pearcy,
`The iterated Aluthge transform of an operator,'
{\em Integr. Equ. Oper. Theory} {\bf 45} (2003), 375-387.

\bibitem{HS2} U.\ Haagerup and H.\ Schultz, 'Invariant subspaces for operators in a general II$_1$--factor',
preprint, arXiv:math.OA/0611256.

\bibitem{HP} D.L.\ Hanson and G.\ Pledger,
'On the mean ergodic theorem for weighted averages',
{\em Z. Wahrscheinlichkeitstheorie verw. Geb.}\ {\bf 13} (1969), 141-149.

\bibitem{P83} K.\ Petersen,
{\em Ergodic Theory},
Cambridge University Press, Cambridge, 1983.

\bibitem{W03} D.\ Wang,
`Heinz and McIntosh inequalities, Aluthge transormation and the spectral measure,'
{\em Math.\ Inequal.\ Appl.}\ {\bf 6} (2003), 121-124.

\bibitem{Y02} T.\ Yamazaki,
`An expression of spectral radius via Aluthge transformation,'
{\em Proc.\ Amer.\ Math.\ Soc.}\   {\bf 130} (2002), 1131-1137.

\end{thebibliography}

\end{document}